\documentclass{article}
\usepackage{amsmath, amsfonts, enumerate,amsthm, cleveref}
\usepackage{graphicx}
\usepackage{caption}
\usepackage{subcaption}

\newtheorem{theorem}{Theorem}
\newtheorem{lemma}{Lemma}
\newtheorem{corollary}{Corollary}
\newtheorem{definition}{Definition}
\newtheorem{proposition}{Proposition}

\title{A geometrical approach to measure irrationality}
\author{Pedro Morales-Almazan\\ \texttt{pmorales@math.utexas.edu}}

\begin{document}

\maketitle
\begin{abstract}
We present a geometric way of describing the irrationality of a number using the area of a circular sector $A(r)$. We establish a connection between this and the continued fraction expansion of the number, and prove bounds for $A(r)$ as $r\to\infty$ by describing the asymptotic behavior of the ratios of the denominators of the convergents.
\end{abstract}

\section{Introduction}
The idea of measuring how well irrationals behave under rational approximations is an old subject that have caught the attention of many people over history. 

Since rationals are dense in the reals, we have that every real number can be approximated through sequences of rationals. The idea of measuring the level of irrationality of a number is to describe how well these approximation can work whenever we try to use the smallest numbers possible. Dirichlet showed that we can always find rationals within a close range of a real number $\alpha$ \cite{miller2006invitation,sally2007roots},
\begin{theorem}[Dirichlet]\label{thm:dirichlet}
If $\alpha$ is irrational, then there are infinitely many rational numbers $p/q$ such that 
\begin{equation}
\Big|\alpha-\frac{p}{q}\Big|<\frac{1}{q^2}\,.
\end{equation}
\end{theorem}
Describing how close can one have rationals controlling the size of the denominator has given a way to understand how well can we approximate irrationals using rational numbers. It turns out to be that the exponent 2 in \Cref{thm:dirichlet} is a threshold for algebraic numbers of degree 2, meaning that any exponent $2+\epsilon$, $\epsilon>0$ will make the inequality satisfied only by a finite number of rationals \cite{sally2007roots}.

Liouville later showed that something similar happened with algebraic numbers of any degree. He stated that for an algebraic number $\alpha$ of degree $n$, there is a constant $C(\alpha)$ such that we can never get closer than $C(\alpha)/q^n$ for any rational $p/q$ \cite{apostol1997modular}. This result led Liouville to investigate about very badly approximable numbers, later called \emph{Liouville numbers} \cite{Liouville51,sally2007roots}.

Following this idea of measuring how well can one approximate an irrational based on the number of solutions to the inequalities described before, it is customary to define the \emph{irrationality measure} or \emph{irrationality exponent} of a real number $\alpha$ to be the smallest real number $\mu$ such that for every $\epsilon > \mu$,
\begin{equation}
\Big|\alpha-\frac{p}{q}\Big|\leq \frac{1}{q^\epsilon}
\end{equation}
has at most a finite number of rational solutions $p/q$ \cite{sally2007roots, Sondow04,MTK:6647164,hardy2008introduction}. Roth established that $\mu$ is in fact 2 for all algrebraic numbers of any degree \cite{MTK:6647164,sally2007roots}. 

In this paper we propose a geometric way of measuring irrationality that captures how well can we approximate real numbers using rationals. In \Cref{sec:geo} we give the geometric interpretation of the approximation problem for a real number $\alpha$ as the area $A(r)$ of a circular sector determined by integer lattice points, and make the connection with the continued fraction expansion of $\alpha$. In \Cref{sec:bounds} we exploit the properties of the convergents for a real number together with the geometric properties derived from the interpretation developed in \Cref{sec:geo} to get general bounds for the local maxima $\{M_k\}$ and minima $\{m_k\}$ of $A(r)$. \Cref{sec:alg2} gives results about the sequence of ratios for the denominators of convergents for algebraic numbers of degree 2 and give bounds for $\{m_k\}$ and $\{M_k\}$ at infinity. Finally, in \Cref{sec:gen} we describe the behavior of general numbers and state a characterization of numbers with unbounded coefficients in its continued fraction expansion in terms of the behavior of $A(r)$.

\section{From geometry to continued fractions}\label{sec:geo}
To investigate the irrationality of a number $\alpha\in\mathbb{R}^+$, consider the the line $L$ given by $y=\alpha x$ in $\mathbb{R}^2$. For a given $r\in\mathbb{R}^+$ consider $S_r$ to be the biggest circular sector of radius $r$ centered at the origin and symmetric around $L$ that does not contain any point from the integer lattice $\mathbb{Z}^2$. The idea is to investigate the area $A(r)$ of $S_r$ as $r\to\infty$.

\begin{figure}[h]
\centering
\includegraphics[height=2in]{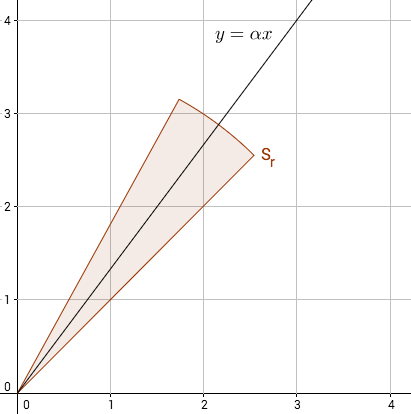}
\caption{Circular sector $S_r$}
\end{figure}

Let $\theta(r)$ be the aperture angle of $S_r$. That is, 
\begin{equation}
A(r)=\frac{r^2}{2}\theta(r)\,.
\end{equation}
Notice that $\theta(r)=2\pi$ for small $r$. When $r$ starts increasing, $\theta(r)$ decreases and we have that there is at least a point $(q,p)\in\mathbb{Z}^2$ on the border of the circular sector. Without loss of generality, suppose that the point lies on either of the straight sides of $S_r$ and not on the arc. This will make $\theta(r)$, and hence $A(r)$, to be left continuous.

The point $(q,p)$ is then located in the first quadrant and is such that $(q,p)$ gives the best rational approximation to the slope of $L$ inside the circle of radius $r$ with center at the origin. Thus we have that 
\begin{equation}\label{eq:angular}
\theta(r)=2\big|\arctan(p/q)-\arctan(\alpha)\big|\,.
\end{equation}
Notice that if there is another such integer point $(q',p')$ for $p'^2+q'^2<r^2$, it is either completely outside of $S_r$ or it lies on either of the straight sides of $S_r$, and we have that \cref{eq:angular} is well defined, as $\arctan(p/q)-\arctan(\alpha)=\pm(\arctan(p'/q')-\arctan(\alpha))$.

We also have that $(q,p)$ provides the best rational approximation of $\alpha$ with \emph{small} numbers. That is, we have the following result \cite{10.2307/2324717,brasselet2007singularities,karpenkov2013geometry},

\begin{proposition}
Let $r>0$ and $(q,p)$ as above. Then $p/q$ is the closest rational to $\alpha$ in the set
\begin{equation}
B_r=\{(q,p):p^2+q^2<r^2\}\,.
\end{equation}
\begin{proof}
Suppose that there is a better approximation to $\alpha$ in $B_r$. That is, let $(q',p')\in B_r$ such that 
\begin{equation}
\Big|\alpha-\frac{p}{q}\Big|>\Big|\alpha-\frac{p'}{q'}\Big|\,.
\end{equation}
Then, we have the line $L'$ given by $\displaystyle y=\frac{p'}{q'}x$ intersects the interior of $S_r$ since 
\begin{equation}
\Big|\alpha x-\frac{p}{q}x\Big|>\Big|\alpha x-\frac{p'}{q'}x\Big|\,,
\end{equation}
and hence $(q',p')$ is in the interior of $S_r$ which is a contradiction.
\end{proof}
\end{proposition}

This result suggest that in order to analyze the behavior of $A(r)$ we need to study the best approximations to $\alpha$ by rationals in $B_r$. This can be achieved by following the continued fraction expansion of $\alpha$, which gives the best rational approximations by means of its convergents \cite{lange1999continued, jones2009continued, hardy2008introduction, stark1978introduction}.

\begin{lemma}\label{lem:best}
Let $\alpha=[a_0;a_1,a_2,\dots]$ be the simple continued fraction expansion of $\alpha$. Then the convergent $h_k=p_k/q_k=[a_0;a_1,a_2,\dots,a_k]$, $\gcd(pk,q_k)=1$, is the best rational approximation to $\alpha$ with denominator smaller or equal than $q_k$.
\end{lemma}

\Cref{lem:best} implies that $\theta$ can be piece-wise defined as
\begin{equation}
\theta(r)=
\begin{cases}
2\pi\,, & 0<r\leq r_0\\
2\Big|\arctan(h_k)-\arctan(\alpha)\Big|\,, & r_{k}< r\leq r_{k+1}\,,
\end{cases}
\end{equation}
where we define the radii $r_k$ in terms of the convergents for $\alpha$,
\begin{equation}
r_{k}=\sqrt{p_k^2+q_k^2}\,.
\end{equation}

With this, we can compute the area  $A(r)$ of the circular sector using the convergents of $\alpha$.

\begin{proposition}\label{pro:fun}
We have that for $r>0$, the area of the circular sector $S_r$ is given by
\begin{equation}
A(r)=
\begin{cases}
\pi r^2\,, & 0<r\leq r_0\\
r^2\Big|\arctan(h_k)-\arctan(\alpha)\Big|\,, & r_{k}< r\leq r_{k+1}\,.
\end{cases}
\end{equation}
\end{proposition}

\begin{figure}[ht!]
    \centering
    \begin{subfigure}[t]{0.5\textwidth}
        \centering
        \includegraphics[height=1.2in]{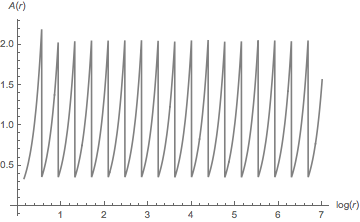}
        \caption{$\alpha=\sqrt{2}$}
    \end{subfigure}%
    ~ 
    \begin{subfigure}[t]{0.5\textwidth}
        \centering
        \includegraphics[height=1.2in]{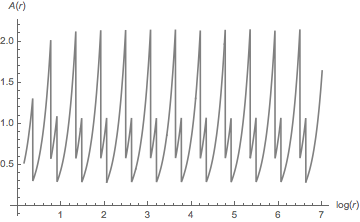}
        \caption{$\alpha=\sqrt{3}$}
    \end{subfigure}\\
    \centering
    \begin{subfigure}[t]{0.5\textwidth}
        \centering
        \includegraphics[height=1.2in]{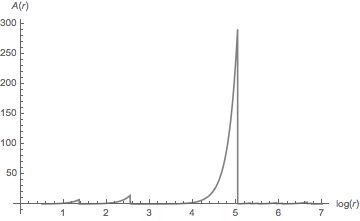}
        \caption{$\alpha=\pi$}
    \end{subfigure}%
    ~ 
    \begin{subfigure}[t]{0.5\textwidth}
        \centering
        \includegraphics[height=1.2in]{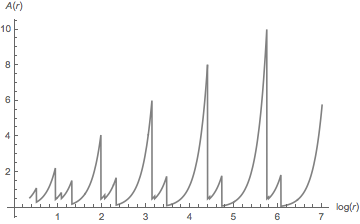}
        \caption{$\alpha=e$}
    \end{subfigure}\\   
    \caption{Graph of $A(r)$ for different values of $\alpha$.}
    \label{fig:areas}
\end{figure}

Some graphs of $A(r)$ for different $\alpha$ are shown in \Cref{fig:areas}.

\section{Bounds for $A(r)$}\label{sec:bounds}
For $\alpha\in\mathbb{Q}^+$ we have that $A(r)$ eventually vanishes. Hence we will focus in the case of irrationals $\alpha\in\mathbb{R}^+$. Our goal now is to describe the behavior of $A(r)$ in the limit $r\to\infty$. This is non trivial, as the quadratic growth of $r$ could compensate for the decrease in the aperture angle $\theta$. We show in this section that the area remains bounded and provide estimates for the bounds.

By \Cref{pro:fun}, we have that $A(r)$ is piece-wise increasing and satisfies
\begin{multline}
m_k=r_k^2\Big|\arctan\left(h_k\right)-\arctan(\alpha)\Big|<A(r)\\
\leq M_k=r_{k+1}^2\Big|\arctan\left(h_k\right)-\arctan(\alpha)\Big|\,,
\end{multline}
for $r_k<r\leq r_{k+1}$ and $k\geq0$.

In order to get bounds for $A(r)$, we use bounds on the error given by the convergents of $\alpha$. For simple continued fractions we have the following result \cite{stark1978introduction, hardy2008introduction, wall1948analytic}
\begin{lemma}\label{lem:conv}
Let $\{h_k\}$ be the convergents for $\alpha$. Then 
\begin{equation}
h_0<h_2<\dots<\alpha<\dots<h_3<h_1\,,
\end{equation}
and the differences between convergents are given by
\begin{equation}
h_{k+1}-h_{k}=\frac{(-1)^k}{q_kq_{k+1}}\,,\qquad h_{k+2}-h_{k}=\frac{(-1)^ka_{k+2}}{q_kq_{k+2}}\,.
\end{equation}
\end{lemma}
The structure of odd and even convergents together with the difference between consecutive, consecutive even, and consecutive odd convergents allow us to get lower and upper estimates for the approximation error given by convergents \cite{stark1978introduction, hardy2008introduction, wall1948analytic}.
\begin{lemma}\label{lem:apboun}
\begin{equation}\label{eq:aboun}
\frac{a_{k+2}}{q_{k+2}q_k}<\Big|\alpha-h_k\Big|<\frac{1}{q_kq_{k+1}}\,.
\end{equation}
\begin{proof}
Since $\alpha$ is between $h_k$ and $h_{k+1}$, we have that the upper bound follows. Similarly, since $|\alpha-h_k|<|\alpha-h_{k+2}|$ and both $h_k,h_{k+2}$ on the same side of $\alpha$, we have that $|\alpha-h_k|$ is smaller than $|h_{k+2}-h_k|$ and the lower bound follows.
\end{proof}
\end{lemma}
With \Cref{lem:apboun} we can show a couple of inequalities that are important in order to find bounds for the area of the circular sector.
\begin{theorem}\label{thm:epboun}
For $\alpha=[a_0;a_1,\dots]$ with convergents $h_k=p_k/q_k$ and convergent radii $r_k^2=p_k^2+q_k^2$, we have that
\begin{equation}
\left(h_k^2+1\right)\left(\frac{q_k}{q_{k+1}}\right)\left(\frac{1}{1+\frac{1}{a_{k+2}}\left(\frac{q_k}{q_{k+1}}\right)}\right)<r_k^2\Big|\alpha-h_k\Big|
\end{equation}
and
\begin{equation}
r_k^2\Big|\alpha-h_k\Big|<\left(h_k^2+1\right)\left(\frac{q_k}{q_{k+1}}\right)
\end{equation}
\begin{proof}
Using \Cref{lem:apboun} we get the upper bound by multiplying \cref{eq:aboun} by $r_k^2=p_k^2+q_k^2$. For the lower bound, it is known in the literature \cite{stark1978introduction, hardy2008introduction, wall1948analytic} that $\{p_k\}$ and $\{q_k\}$ satisfy recurrence relations 
\begin{equation}\label{eq:rec}
p_k=a_kp_{k-1}+p_{k-2}\,,\qquad q_k=a_kq_{k-1}+q_{k-2}\,,
\end{equation}
with initial conditions
\begin{equation}
p_0=a_0\,,\qquad p_1=a_0a_1+1\,,
\end{equation}
and 
\begin{equation}
q_0=1\,,\qquad q_1=a_1\,.
\end{equation}
Therefore using the recurrence relation for $q_{k+2}$ and substituting in \cref{eq:aboun} gives the desired lower bound.
\end{proof}
\end{theorem}
This result bring us closer to find bounds for $m_k$, we only need a technical lemma before we can state these.
\begin{lemma}\label{lem:arctan}
For $x>0$ and $0<\epsilon\ll 1$,
\begin{equation}\label{eq:arctanineq}
\frac{\epsilon}{1+x^2}-\frac{3\sqrt{3}}{16}\epsilon^2<\Big|\arctan(x\pm\epsilon)-\arctan(x)\Big|<\frac{\epsilon}{1+x^2}+\frac{3\sqrt{3}}{16}\epsilon^2\,.
\end{equation}
\begin{proof}
Fix $x>0$ and let $0<\epsilon\ll1$ such that the Maclaurin expansion for $\arctan(x+\epsilon)$ as a function of $\epsilon$ converges. 
Thus, we have that 
\begin{equation}
\arctan(x+\epsilon)=\arctan(x)+\frac{\epsilon}{1+x^2}+R_2(x,\epsilon)\,,
\end{equation}
where $R_2$ is the truncation error. By Taylor's theorem and the Mean value theorem, we have that 
\begin{equation}
R_2(x,\epsilon)=-\frac{\epsilon^2(x+y)}{(1+(x+y)^2)^2}\,,
\end{equation}
for some $y\in(0,\epsilon)$. Since $R_2(x,\epsilon)$ has global extrema when $x+y=\pm\frac{1}{\sqrt{3}}$, we have that 
\begin{equation}
\Big|R_2(x,\epsilon)\Big|\leq \frac{3\sqrt{3}}{16}\,.
\end{equation}
With this, we have that 
\begin{equation}
\arctan(x+\epsilon)-\arctan(x)=\frac{\epsilon}{1+x^2}-\frac{\epsilon^2(x+y)}{(1+(x+y)^2)^2}\,, 
\end{equation}
and we have the bounds
\begin{equation}
\frac{\epsilon}{1+x^2}-\frac{3\sqrt{3}}{16}\epsilon^2\leq\arctan(x+\epsilon)-\arctan(x)\leq\frac{\epsilon}{1+x^2}\,.
\end{equation}
Similarly, we have that 
\begin{equation}
\arctan(x)-\arctan(x-\epsilon)=\frac{\epsilon}{1+x^2}+\frac{\epsilon^2(x+y')}{(1+(x+y')^2)^2}\,,
\end{equation}
for some $y'\in(0,\epsilon)$ and thus we have that
\begin{equation}
\frac{\epsilon}{1+x^2}\leq\arctan(x)-\arctan(x-\epsilon)\leq\frac{\epsilon}{1+x^2}+\frac{3\sqrt{3}}{16}\epsilon^2\,
\end{equation}
and the result follows.
\end{proof}
\end{lemma}
Hence, using \Cref{thm:epboun,lem:arctan} we are able to bound $m_k$ from above and below,
\begin{theorem}\label{thm:boun}
Let $\epsilon_k=\alpha-h_k$. Then
\begin{equation}
m_k<\left(\frac{h_k^2+1}{\alpha^2+1}\right)\left(\frac{q_k}{q_{k+1}}\right)\left(1+\frac{3\sqrt{3}}{16}(\alpha^2+1)\epsilon_k\right)\,,
\end{equation}
and
\begin{equation}
\left(\frac{h_k^2+1}{\alpha^2+1}\right)\left(\frac{q_k}{q_{k+1}}\right)\left(\frac{1}{1+\frac{1}{a_{k+2}}\left(\frac{q_k}{q_k+1}\right)}-\frac{3\sqrt{3}}{16}(\alpha^2+1)\epsilon_k\right)<m_k\,.
\end{equation}
\begin{proof}
Setting $\epsilon=\epsilon_k$ in \Cref{lem:arctan}, multiplying \cref{eq:arctanineq} by $r_k^2$, and using \Cref{thm:epboun} together with the fact that $\epsilon_k>0$ lead to the result.
\end{proof}
\end{theorem}
This result gives bounds for $m_k$ for each $k$, in terms of the convergents of $\alpha$. To study the behavior of $m_k$ as $k\to\infty$ we need a lemma involving bounds for sequences.
\begin{lemma}\label{lem:lim}
Let $x_n<y_n$ be sequences of real numbers such that $y_n\to y$. Then 
\begin{equation}
\limsup_{n\to\infty}x_n\leq y.
\end{equation}
\begin{proof}
Suppose by contradiction that 
\begin{equation}
y<\limsup_{n\to\infty} x_n\,.
\end{equation}
Then, there is a sequence os positive integers $\{n_k\}$ such that 
\begin{equation}
y<\lim_{k\to\infty} x_{n_{k}}\,,
\end{equation}
but we have that $x_{n_{k}}<y_{n_k}$ and hence
\begin{equation}
y<\lim_{k\to\infty} x_{n_k}<\lim_{k\to\infty}y_{n_{k}}=y\,,
\end{equation}
which is a contradiction.
\end{proof}
\end{lemma}
Therefore, we have that \Cref{thm:boun} together with the fact that $\epsilon_k\to 0$ as $k\to\infty$ gives that we can bound the behavior of $m_k$ at infinity,
\begin{corollary}
\begin{equation}
\limsup_{k\to\infty} m_k\leq 1\,.
\end{equation}
\begin{proof}
We have that $\{q_k\}$ is an increasing sequence so that $q_{k-1}/q_k<1$. Using this and \Cref{lem:lim} gives the result.
\end{proof}
\end{corollary}

Similarly, we can find bounds for $M_k$ using  \Cref{lem:apboun,lem:arctan}. Using the same approach as before, we get the following result.

\begin{theorem}\label{thm:M}
Using the notation as before, 
\begin{equation}
M_k<\left(\frac{h_{k+1}^2+1}{\alpha^2+1}\right)\left(\frac{q_{k+1}}{q_k}\right)\left(1+\frac{3\sqrt{3}}{16}(\alpha^2+1)\epsilon_k\right)\,,
\end{equation}
and
\begin{equation}
\left(\frac{h_{k+1}^2+1}{\alpha^2+1}\right)\left(\frac{q_{k+1}}{q_k}\right)\left(\frac{1}{1+\frac{1}{a_{k+2}}\left(\frac{q_k}{q_{k+1}}\right)}-\frac{3\sqrt{3}}{16}(\alpha^2+1)\epsilon_k\right)<M_k\,.
\end{equation}
\begin{proof}
Multiplying \cref{eq:aboun} by $r_{k+1}^2=p_{k+1}^2+q_{k+1}^2$, using \Cref{lem:conv}, and setting $\epsilon=\epsilon_{k}=|\alpha-h_{k}|$ in \Cref{lem:arctan} give the result.
\end{proof}
\end{theorem} 

\section{Asypmtotic periodicity of extreme values of $A(r)$ for algebraics of degree 2}\label{sec:alg2}
In this section we analyze the behavior of $A(r)$ for $\alpha$ and algebraic number of order 2. We proceed by completely describing the behavior of $m_k$ and $M_k$ as $k\to\infty$.

In order to start, we use a very useful characterization of these numbers given by Lagrange \cite{davenport1999higher, khinchin1964continued, karpenkov2013geometry}.

\begin{theorem}
If $\alpha$ is an algebraic number of order 2, its continued fraction expansion is eventually periodic.
\end{theorem}

With this result we are able to write any such $\alpha$ as the simple continued fraction
\begin{equation}\label{eq:alg2}
\alpha=[a_0;a_1,a_2,\dots, a_l, \overline{b_1,\dots, b_n}]\,,
\end{equation}
were the bar indicates the repeating part. 

Notice that the bounds given in \Cref{thm:boun,thm:M}, for $m_k$ and $M_k$ respectively, depend on the consecutive ratios of the denominators of the convergents for $\alpha$. Thus we need to describe how these ratios behave and their limit properties.

Before we describe this behavior, we need a couple of technical lemmas over simple continued fractions over the complex numbers,
\begin{equation}\label{eq:cont_frac}
\cfrac{1}{b_1+\cfrac{1}{b_2+\cfrac{1}{b_3+\ddots}}}\,,
\end{equation}
where $b_k\in\mathbb{C}$. The first result gives a sufficient condition for \cref{eq:cont_frac} in order to get convergence, while the second one states how sensitive is the convergence of the continued fraction depending on the coefficients $b_n$.

We start by defining the idea of uniform convergence region and simple region \cite{lange1999continued, jones2009continued, Lange1999355}. 

\begin{definition}[Convergence regions]
Let $f_k$ be the $k$th convergent of a continuous fraction given by \cref{eq:cont_frac}, and let $\lim_{k\to\infty}f_k=f$ whenever the limit exists. 

Let $\{\Omega_k\}$ be a sequence of regions in the complex plane. If $\{\Omega_k\}$ is such that if $b_k\in\Omega_k$ for all $k$ then the continued fraction \cref{eq:cont_frac} converges, we say that $\{\Omega_k\}$ is a convergence sequence for the continued fraction \cref{eq:cont_frac}.

We say that the convergence sequence $\{\Omega_k\}$ is uniform if there are positive real numbers $\epsilon_k$ depending only on the sequence, such that $\lim_{k\to\infty}\epsilon_k=0$ and $|f-f_k|<\epsilon_k$ whenever $b_k\in\Omega_k$.
\end{definition}

\begin{definition}[Uniform simple convergence region]
If $\{\Omega_k\}$ is a Uniform convergence sequence such that $\Omega_k=\Omega$ for all $k$, we call this a uniform simple convergence region.
\end{definition}

The first of these results is due to Van Vleck, and gives a sufficient condition to ensure convergence of \cref{eq:cont_frac} \cite{10.2307/1986206}.

\begin{lemma}[Van Vleck]\label{lem_vanvleck}
If the series $\sum|b_p|$ diverges, then the continued fraction
\begin{equation}
\cfrac{1}{b_1+\cfrac{1}{b_2+\cfrac{1}{b_3+\ddots}}}
\end{equation}
converges.
\end{lemma}

The second result deals gives regions in which the convergence \cref{eq:cont_frac} is uniform on its coefficients \cite{lange1999continued}.

\begin{lemma}[Uniform Circle Theorem]\label{thm:uniform}
Let $c$ be a real number, $r=\sqrt{1+c^2}$, and $B[c]:=\{z:|z+2c|\geq 2r\}$. Then $B[c]$ is a best uniform simple convergence region for  \cref{eq:cont_frac}.
\end{lemma}

With these results, we are in good position to prove the following result which will help us later in order to find the behavior of $m_k$ and $M_k$ when $k\to\infty$.

\begin{theorem}\label{thm:cont}
Let $\{a_{i,n}\}_{n,i=0}^\infty$ be sequences of real numbers greater or equal than 1 such that 
\begin{equation}
\lim_{n\to\infty} a_{i,n}=a_i\,,
\end{equation}
for all $i\in\mathbb{N}$.
If $b_n=[a_{0,n};a_{1,n},a_{2,n},\dots]$, then
\begin{equation}
\lim_{n\to\infty} b_n=b=[a_0;a_1,a_2,\dots]\,.
\end{equation}

\begin{proof}
For convenience, denote a finite continued fraction $[c_0;c_1,c_2,\dots,c_l]$ by $K_{n=0}^l(1/c_n)$ and analogously $[c_0;c_1,c_2,\dots]$ by $K_{n=0}^\infty(1/c_n)$.

By \Cref{lem_vanvleck}, we have that $b=K_{n=0}^\infty(1/a_n)$ converges since $\sum|a_n|$ diverges. Now, to prove that the limit exists, let $\epsilon>0$. By \Cref{thm:uniform}, there exists $N\in\mathbb{N}$ such that $0<\epsilon_k<\epsilon/3$ for all $k>N$. Hence for all $n$, 
\begin{equation}
\Big|b_n-K_{i=0}^k(1/a_{i,n})\big|<\epsilon_k<\epsilon/3\,.
\end{equation}

Notice that $K_{i=0}^k(1/x_i)$ is a rational function on $x_i$ and it is uniformly continuous for compact subsets away from the singularities. Since $a_i>0$ we have that $x_i=a_i$ are not singularities and hence there exists $M$ such that 
\begin{equation}
\Big|K_{i=0}^k(1/a_{i,n})-K_{i=0}^k(1/a_{i,m})\Big|<\epsilon/3\,,
\end{equation}
for all $n,m>M$.

Therefore, we have that 
\begin{multline}
\Big|b_n-b_m\Big|\leq \Big|b_n-K_{i=0}^k(1/a_{i,n})\Big|+\Big|b_m-K_{i=0}^k(1/a_{i,m})\Big|\\
+\Big|K_{i=0}^k(1/a_{i,m})-K_{i=0}^k(1/a_{i,n})\Big|<\epsilon\,,
\end{multline}
for all $n,m>M$. Then, we have that $\{b_k\}$ is a Cauchy sequence and hence the limit exists and it is equal to $b$.
\end{proof}
\end{theorem}

This result says that a \cref{eq:cont_frac} when viewed as a function on the variables $b_i$ is continuous in each variable. With this, we are able to describe the behavior of the partial quotients for algebraic numbers of order 2. 

Let $c_k$ be the ratio of consecutive partial quotients,
\begin{equation}
c_k=\frac{q_k}{q_{k-1}}\,.
\end{equation}

\begin{theorem}\label{thm:periods}
If $\alpha$ has a periodic continued fraction expansion of length $n$, then $\{c_k\}$ is asymptotically periodic with $n$ different subsequential limits.

\begin{proof}
Using the recurrence relation \cref{eq:rec}, we have that 
\begin{equation}
c_k=\frac{q_k}{q_{k-1}}=a_k+\frac{q_{k-2}}{q_{k-1}}=a_k+\frac{1}{c_{k-1}}\,.
\end{equation}
Therefore, by an inductive argument, we have that 
\begin{equation}
c_k=[a_k;a_{k-1},a_{k-2},\dots, a_0]\,.
\end{equation}
Let
\begin{equation}
\alpha=[a_0;a_1,a_2,\dots, a_l, \overline{b_1,b_2,b_3,\dots,b_n}]\,,
\end{equation}
then we have that 
\begin{equation}
c_{l+kn+i}=[b_i;b_{i-1},b_{i-2},\dots,b_1,\underbrace{b_n,b_{n-1},\dots,b_{2},b_1}_{k \text{ times}},a_l,a_{l-1},\dots, a_1,a_0]\,.
\end{equation} 
Thus, by \Cref{thm:cont} we have that 
\begin{equation}
C_i=\lim_{k\to\infty}c_{l+kn+i}=[\overline{b_i;b_{i-1},b_{i-2},\dots,b_{i+2},b_{i+1}}]\,,
\end{equation}
for $1< i\leq n$.
\end{proof}
\end{theorem}

With this and \Cref{thm:boun}, we get that 
\begin{theorem}
\begin{equation}
\frac{1}{C_i+\frac{1}{b_{i+1}}}\leq \liminf_{k\to\infty}m_{l+kn+i-1}\leq \frac{1}{C_i}\,,
\end{equation}
with $1\leq i\leq n$, $C_i$ as in \Cref{thm:periods}, and $b_{n+1}=b_1$.
\end{theorem}
An immediate result can be obtained from here by taking the extremes of all such bounds to get an overall bound for the sequence of infima $m_k$.
\begin{corollary}
\begin{equation}
m\leq \liminf_{k\to\infty} m_k\leq\nu \,,
\end{equation}
where $m$ is the minimum 
\begin{equation}
m=\min \left\{\frac{1}{C_i+\frac{1}{b_{i+1}}}:i\in\{1,2,3,\dots,n\}\right\}\,,
\end{equation}
and $\nu$ is the maximum 
\begin{equation}
\nu=\max \left\{\frac{1}{C_i}:i\in\{1,2,3,\dots,n\}\right\}\,.
\end{equation}

\end{corollary}

Likewise, we can bound the behavior of the maxima $M_k$ of the areas $A(r)$. Using the periodicity of algebraics of degree two, we can show that these are also bounded and the bounds depend on the continued fraction expansion of $\alpha$.

\begin{theorem}
With $\alpha$ given by \cref{eq:alg2} and $C_i$ as in \Cref{thm:periods}, we have that the maxima $M_k$ satisfy
\begin{equation}
\limsup_{k\to\infty} M_{l+kn+i-1}\leq C_i\,,
\end{equation}
and
\begin{equation}
\frac{C_i^2}{C_i+\frac{1}{b_{i+1}}}\leq \liminf_{k\to\infty} M_{l+kn+i-1}
\end{equation}
\begin{proof}
Using the values for the subsequential limits found in \Cref{thm:periods} into \Cref{thm:M} together with \Cref{lem:lim} gives the result.
\end{proof}
\end{theorem} 
For this, we can take the minimum and maximum values of the bounds appearing in the previous result to get general bounds to the local maxima of $A(r)$.
\begin{corollary}
Let 
\begin{equation}
M=\max\left\{C_i:i\in\{1,2,3,\dots,n\}\right\}\,,
\end{equation}
and
\begin{equation}
\mu=\min\left\{\frac{C_i^2}{C_{i}+\frac{1}{b_{i+1}}}:i\in\{1,2,3,\dots,n\}\right\}\,.
\end{equation}
Then we have that 
\begin{equation}
\limsup_{k\to\infty} M_k\leq M\,,
\end{equation}
and
\begin{equation}
\mu\leq \liminf_{k\to\infty} M_k\,. 
\end{equation}
\end{corollary}

With these results we are able to bound the behavior at infinity of algebraics of degree 2.

\begin{corollary}
\begin{equation}
\limsup_{r\to\infty}A(r)\leq M\,,
\end{equation}
and
\begin{equation}
m\leq \liminf_{r\to\infty}A(r)\,.
\end{equation}
\end{corollary}

\section{Bounds of $A(r)$ for general numbers}\label{sec:gen}
On analyzing general real numbers, we can establish some results about the relation between the growth of their continued fraction coefficients and behavior of the sequences $\{m_k\}$ and $\{M_k\}$.

First, as described in \Cref{thm:periods}, we can describe the growth of the ratios of denominators $\{c_k\}$ with the growth of the continued fraction coefficients $\{a_k\}$.

\begin{theorem}\label{thm:order}
For irrationals $\alpha\in\mathbb{R}^+$, we have that 
\begin{equation}
c_k= a_k+O(1)
\end{equation}
as $k$ tends to infinity.
\begin{proof}
Since $\{q_k\}$ is an increasing sequence of positive integers, we have that $c_k>1$ and therefore
\begin{equation}
c_k-a_k=\frac{q_k}{q_k-1}-a_k=\frac{q_{k-2}}{q_{k-1}}<1\,,
\end{equation}
for all $k$. Thus the result follows.
\end{proof}
\end{theorem}

With this result, we immediately have a very interesting result about the sequences of extrema $\{m_k\}$ and $\{M_k\}$ for $\alpha$.

\begin{theorem}
Let $\alpha=K_{k=0}^\infty(1/a_k)$ be a positive irrational number. Then, we have that 0 is a subsequential limit of $\{m_k\}$ and $\{M_k\}$ is unbounded if and only if $\{a_k\}$ is unbounded.

\begin{proof}
If $\{a_k\}$ is unbounded, we have that by \Cref{thm:order} 
\begin{equation}
\limsup_{k\to\infty}c_k=\infty\,.
\end{equation}
This and \Cref{thm:boun} give that 
\begin{equation}
\liminf_{k\to\infty} m_k=0\,,
\end{equation}
and by \Cref{thm:M} we have that 
\begin{equation}
\limsup_{k\to\infty} M_k=\infty\,,
\end{equation}
and the sufficient part follows.
To prove necessity, by contraposition, suppose that $\{a_k\}$ is bounded. Then we have that $\{c_k\}$ is bounded and \Cref{thm:boun,thm:M} show that 
\begin{equation}
\liminf_{k\to\infty} m_k>0
\end{equation}
and
\begin{equation}
\limsup M_k<\infty\,,
\end{equation}
and this completes the proof.
\end{proof}
\end{theorem}

This result says that if $\alpha$ has unbounded continued fraction coefficients, then $A(r)$ oscillates between 0 and $\infty$ as $r\to\infty$. For example, since the continued fraction expansion for e is given by \cite{10.2307/2318113}
\begin{equation}
e=[2;1,2,1,1,4,1,1,6,1,1,8,1,1,\dots]\,,
\end{equation}
we have that for any $a>0$ and for any $R>0$ there are $r_m,r_M>R$ such that 
\begin{equation}
A(r_m)<a\quad \text{ and }\quad A(r_M)>a\,.
\end{equation}

\bibliographystyle{plain}
\bibliography{references}

\end{document}